\documentclass[11pt]{amsart}
\usepackage{amsmath,amstext,amsthm,amsfonts,amssymb,psfrag}

\theoremstyle{plain}

\newtheorem{theorem}{Theorem}[section]
\newtheorem{lemma}[theorem]{Lemma}
\newtheorem{corollary}[theorem]{Corollary}
\newtheorem{proposition}[theorem]{Proposition}
\theoremstyle{remark}

\newtheorem{remark}[theorem]{Remark}

\title{Topological rigidity for holomorphic foliations}
\author{MAHDI TEYMURI GARAKANI}

\address{Max-Planck-Institut fur Mathematik, Vivatsgasse 7 Bonn
D-53111 Germany}

\email{teymuri@impa.br}

\keywords{Singular holomorphic foliation, deformation, unfolding,
rigidity, holonomy }

\subjclass[2000]{Primary 37F75, 58k45}

\begin{document}

\maketitle

\begin{abstract}
We study analytic deformations and unfoldings of holomorphic
foliations in complex projective plane $\mathbb{C}P(2)$. Let
$\{\mathcal{F}_t\}_{t \in \mathbb{D}_{\varepsilon}}$ be
topological trivial (in $\mathbb{C}^2$) analytic deformation of a
foliation $\mathcal{F}_0$ on $\mathbb{C}^2$. We show that under
some dynamical restriction on $\mathcal{F}_0$, we have two
possibilities: $\mathcal{F}_0$ is a Darboux (logarithmic)
foliation, or $\{\mathcal{F}_t\}_{t \in \mathbb{D}_{\varepsilon}}$
is an unfolding. We obtain in this way a link between the
analytical classification of the unfolding and the one of its
germs at the singularities on the infinity line. Also we prove
that a finitely generated subgroup of
$\mathrm{Diff}(\mathbb{C}^n,0)$ with polynomial growth is
solvable.
\end{abstract}

\section{Introduction}

Let $\mathrm{Fol}(M)$ denote the set of holomorphic foliations on
a complex manifold $M$. An \textit{analytic deformation} of
$\mathcal{F} \in \mathrm{Fol}(M)$ is an analytic family
$\{\mathcal{F}_t\}_{t \in Y}$ of foliations on $M$, with
parameters on an analytic space $Y$, such that there exists a
point $"0" \in Y$ with $\mathcal{F}_0=\mathcal{F}$. Here we will
only consider deformations where $Y=\mathbb{D} \subset \mathbb{C}$
is a unitary disk. A \textit{topological equivalence} (resp.
\textit{analytical equivalence}) between two foliations
$\mathcal{F}_1$ and $\mathcal{F}_2$ is a homeomorphism (resp.
biholomorphism) $\phi :M \rightarrow M$, which takes leaves of
$\mathcal{F}_1$ onto leaves of $\mathcal{F}_2$, and such that
$\phi(\mathrm{Sing}(\mathcal{F}_1)) = \mathrm{Sing}(\mathcal{F}_2)
$. The deformation $ \{{\mathcal{F}_t}\}_{t \in \mathbb{D}} $ is
\textit{topologically trivial} (resp. \textit{analytically
trivial}) if there exists a continuous map (resp. holomorphic map)
$\phi: M \times \mathbb{D} \rightarrow M$, such that each map
$\phi_t = \phi(.,t): M \rightarrow M$ is a topological equivalence
(resp. analytical equivalence) between $\mathcal{F}_t$ and
$\mathcal{F}_0$.

Let $\mathcal{C} \subset Fol(M)$ be a class of foliations. A
foliation $\mathcal{F}_0 \in \mathcal{C}$ is topologically rigid
in the class if any topologically trivial deformation
$\{\mathcal{F}_t\}_{t \in \mathbb{D}}$ of $ \mathcal{F}_{0}$ with
$ \mathcal{F}_t \in \mathcal{C} $ is analytically trivial.

We also say that $ \mathcal{F}_{0} \in \mathcal{C} $ is
$\mathcal{U}$-topological rigid in the class $\mathcal{C} $, where
$\mathcal{U} \subset M $ is an open subset, if any analytic
deformation $\{{\mathcal{F}_t}\}_{t \in \mathbb{D}}$ of $
\mathcal{F}_{0} $  with $\mathcal{F}_{t} \in \mathcal{C}$,
$\forall{t}$; which is topologically trivial in $\mathcal{U}$, is
in fact analytically trivial in $M$.

In this part we will be concerned with holomorphic foliations in
$\mathbb{C}P(2)$. These foliations are motivated by Hilbert
Sixteenth Problem on the number and position of limit cycles of
polynomial differential equations

$$ \frac{dy}{dx}=\frac{P(x,y)}{Q(x,y)} $$

in the real plane $(x,y) \in \mathbb{R}^2$ where $P$ and $Q$ are
relatively prime polynomials. A major attempt in this line was
started in 1956 by a seminal work of I. Petrovski and E. Landis
[16]. They consider $(\star)$ as a differential equations in the
complex plane $(x,y) \in \mathbb{C}^2$, with $t$ now being a
complex time parameter. The integral curves of the vector field
are now either singular points which correspond to the common
zeros of $P$ and $Q$, or complex curves tangent to the vector
field which are holomorphically immersed in $\mathbb{C}^2$. This
gives rise to a holomorphic foliation by complex curves with a
finite number of singular points. One can easily see that this
foliation extends to the complex projective plane
$\mathbb{C}P(2)$, which is obtained by adding a line at infinity
to the plane $\mathbb{C}^2$. Conversely any holomorphic foliation
by curves on $\mathbb{C}P(2)$ is given in an affine space
$\mathbb{C}^2 \hookrightarrow \mathbb{C}P(2)$ by a polynomial
vector field $X=(P,Q) \in \mathfrak{X}(\mathbb{C}^2)$ with
$\mathrm{gcd}(P,Q)=1.$

Although they didn't solve this problem, however they introduced a
truly novel method in geometric theory of ordinary differential
equations. In 1978, Il'yashenko made a fundamental contribution to
the problem. Following the general idea of Petrovski and Landis,
he studied equations $(\star)$ with complex polynomials $P$ and
$Q$ from a topological standpoint without particular attention to
Hilbert's question.

We fix the line at infinity $ L_\infty =\mathbb{C}P(2) \setminus
\mathbb{C}^2 $ and denote by $\mathcal{X}(n)$ the space of
foliations of degree $n \in \mathbb{N}$ which leave invariant $
L_\infty$. Let us denote by $\mathcal{F}(n)$ the space of degree
$n$ foliations on $\mathbb{C}P(2)$ as introduced in [10]. We are
interested in the following question:

\textit{Under which conditions topologically trivial deformations
of a foliation $\mathcal{F} \in \mathcal{F}(n)$ are analytically
trivial?}

A remarkable result of Y. Ilyashenko states topological rigidity
for a residual set of foliations on $\mathcal{X}(n)$ if $n\ge 2$.

More precisely we have:

\begin{theorem}
$\mathrm{[8]}$ For any $n\ge 2$ there exists a residual subset
$\mathcal{I}(n) \subset \mathcal{X}(n)$ whose foliations are
topologically rigid in the class $\mathcal{X}(n)$.
\end{theorem}

This result has been later improved by A. Lins Neto, P. Sad and B.
Scardua as follows:

\begin{theorem}
$\mathrm{ [11]}$ For each $n \ge 2$, $\mathcal{X}(n)$ contains an
open dense subset $\mathcal{R} \subset \mathcal{X}(n) $ whose
foliations are topologically rigid in the class $\mathcal{X}(n)$
\end{theorem}

We stress the fact that in both theorems above we consider
deformations $\{{\mathcal{F}_{t}}\}_{t \in \mathbb{D}}$ in the
class $\mathcal{X}(n)$, that is, ${\mathcal{F}}_{t}$ leaves
invariant $L_{\infty} ,\forall{t} \in \mathbb{D}$; and we assume
topological triviality in $\mathbb{C}P(2)$. This last hypothesis
is slightly relaxed by requiring topological triviality for the
set of separatrices through the singularities at $L_{\infty} $:

\begin{theorem}$\mathrm{ [11]}$ For any $n\ge 2$, $\mathcal{X}(n)$ contains an
open dense subset $\mathcal{SR}\mathrm{ig}(n)$ whose foliations
are $s$-rigid in the class $\mathcal{X}(n)$.
\end{theorem}

According to [11] a foliation $\mathcal{F}_0 \in \mathcal{X}(n)$
is $s$-rigid if for any deformation $\{\mathcal{F}_t\}_{t \in
\mathbb{D}} \subset \mathcal{X}(n)$ of $\mathcal{F}_{0}$ with the
$s$-triviality property that is: If $S_t \subset \mathbb{C}^2 $
denotes the set of separatrices of $\mathcal{F}_{t} $ which are
transverse to $L_{\infty}$ then there exists a continuous family
of maps $\phi_{t}: S_{0} \rightarrow \mathbb{C}^2 $ such that
$\phi_0$ is the inclusion map and $\phi_t$ is a continuous
injection map from $S_{0}$ to $ \mathbb{C}^2$ with
${\phi_t}(S_0)=S_t$; then $\{\mathcal{F}_t\}$ is analytically
trivial.

\begin{remark}
Topological triviality in $\mathbb{C}^2 $ implies $s$-triviality.
\end{remark}

Let us change now our point of view.

A deformation $\{\mathcal{F}_t\}_{t \in \mathbb{D}}$ of a
foliation $\mathcal{F}_0$ on a manifold $\mathrm{M}$, is an
\emph{unfolding} if there exists an analytic foliation
$\tilde{\mathcal{F}}$ on $M \times \mathbb{D}$ with the property
that: $\tilde{\mathcal{F}}\arrowvert_{M \times \{t\}} \equiv
{\mathcal{F}_t}$, $\forall {t} \in \mathbb{D}$. In other words,
\emph{an unfolding is a deformation which embeds into an analytic
foliation.} The trivial unfolding of $\mathcal{F}$ is given by the
$\mathcal{F}_{t}:=\mathcal{F}$, $\forall t \in \mathbb{D}$ and
$\tilde{\mathcal{F}}$ is the product foliation $\mathcal{F} \times
\mathbb{D}$ in $M \times \mathbb{D}$.

Two unfoldings $\{\mathcal{F}_t\}_{t \in \mathbb{D}} $ and
$\{\mathcal{F}^1_t\}_{t \in \mathbb{D}}$ of $\mathcal{F}$ are
topologically equivalent respectively analytically equivalent if
there exists a continuous respectively analytic map $\phi :{M}
\times \mathbb{D} \rightarrow M$ such that each map $\phi_t:M
\rightarrow M$, $\phi_t(p)=\phi(p,t)$, is a topological
respectively analytical equivalence between $\mathcal{F}_t$ and
$\mathcal{F}^1_t$.

An unfolding $\{\mathcal{F}_t\}_{t\in \mathbb{D}}$ of a foliation
$\mathcal{F}_0$ on $M$ is said to be topologically rigid in the
class $\mathcal{C} \subset \mathcal{F}(n)$ if any analytic
unfolding $\{{\mathcal{F}^1_t}\}_{t\in \mathbb{D}}$ of
$\mathcal{F}$ ($\mathcal{F}^1_t \in \mathcal{C}$, $ \forall t$),
which is topologically equivalent to
$\{{\mathcal{F}_\mathrm{t}}\}_{t \in \mathbb{D}}$, is necessarily
analytically equivalent.

These notions rewrite theorems (1.1) and (1.2) as follows:

\begin{theorem} $\mathrm{ [8]}$ For any $n\ge 2$ there exists a residual subset
$\mathcal{I}(n) \subset \mathcal{X}(n)$ whose foliations are
topologically rigid trivial unfolding in the class
$\mathcal{X}(n)$.
\end{theorem}

\begin{theorem}$\mathrm{ [11]}$ For each $n\ge 2$, $\mathcal{X}(n)$ contains an
open dense subset $\mathcal{R}ig(n) \subset \mathcal{X}(n) $ whose
foliations are topologically rigid trivial unfolding in the class
$\mathcal{X}(n)$.
\end{theorem}

We are now in conditions of stating our main results concerning
topological rigidity. We stress the fact that a priori, our
deformations are allowed to move the line $L_{\infty}$ (Theorems
A, B and C) which is $\mathcal{F}_0$-invariant by hypothesis. We
also state results for deformations which are topological trivial
in $\mathbb{C}^2$ not necessarily in $\mathbb{C}P(2)$. Finally, we
may relax the hypothesis of hyperbolicity for
$\mathrm{Sing}(\mathcal{F}_0) \cap L_{\infty}$ by allowing
quasi-hyperbolic singularities (i.e the exceptional divisor is
invariant and all of singularities are of saddle-type) and obtain
in this way a link between the analytical classification of the
unfolding and the one of its germs at the singularities $p \in
\mathrm{Sing}(\mathcal{F}_0) \cap L_{\infty}$.

Our main results are the following:\\
\\
\textbf{Theorem A.}\textit{ Given $n\ge 2$ there exists an open
dense subset $\mathcal{R}ig(n) \subset \mathcal{X}(n) $ such that
any foliation in $\mathcal{R}\mathrm{ig}(n)$ is
$\mathbb{C}^2$-topological rigid: any deformation
$\{{\mathcal{F}_t}\}_{t\in \mathbb{D}}$ of
$\mathcal{F}={\mathcal{F}_0}$ which is topologically trivial in
$\mathbb{C}^2$ must be analytically trivial in $\mathbb{C}P(2)$
for $t \approx 0$. }
\\
\\
\textbf{Theorem B.}\textit{ Let $\{\mathcal{F}_t\}_{t\in
\mathbb{D}}$ be topological trivial (in $\mathbb{C}^2$) analytic
deformation of a foliation $\mathcal{F}_0$ on $\mathbb{C}^2$ such
that:
\begin{enumerate}
\item {$\mathcal{F}_0$ leaves $L_{\infty}$ invariant,}
\item {$\forall p \in \mathrm{Sing}(\mathcal{F}_0) \cap L_{\infty}$, $p$ is
       a quasi-hyperbolic singularity,}
\item {$\mathcal{F}_0$ has degree $ n\ge 2$ and exhibits at least two
       reduced singularities in $\mathrm{L}_{\infty}$.}
\end{enumerate}
Then we have two possibilities:
\begin{itemize}
\item {$\mathcal{F}$ is a Darboux (logarithmic) foliation, or}
\item {$\{{\mathcal{F}_t}\}_{t\in {\mathbb{D}}_{\varepsilon}}$ is an unfolding.}
\end{itemize}}

\textit{In this last case the unfolding is analytically trivial if
and only if given a singularity $p \in
\mathrm{Sing}(\mathcal{F}_0) \cap L_{\infty}$ the germ of the
unfolding $\{{\mathcal{F}_t}\}_{t\in \mathbb{D}}$ at $p$ is
analytically trivial for $t \approx 0$.}
\\
\\
We can rewrite Theorem (B) as follows:
\\
\\
\textbf{Theorem C.}\textit{ Let $\mathcal{F}_0$ be a foliation on
$\mathbb{C}P(2)$ with the following properties:
\begin{enumerate}
\item {$\mathcal{F}_0$ leaves $L_{\infty}$ invariant,}
\item {$\forall p \in \mathrm{Sing}(\mathcal{F}_0) \cap L_{\infty}$, $p$ is a
       quasi-hyperbolic singularity,}
\item {$\mathrm{Sing}(\mathcal{F}_0) \cap L_{\infty}$ has at least two
       reduced singularities.}
\end{enumerate}}

\textit{Given two topologically equivalent unfoldings
$\{{\mathcal{F}_t}\}_{t\in \mathbb{D}}$ and
$\{{\mathcal{F}^1_t}\}_{t\in \mathbb{D}}$  of $\mathcal{F}_0$ we
have that they are analytically equivalent if and only if the
germs of the unfoldings are analytically equivalent at the
singular points $ p \in \mathrm{Sing}(\mathcal{F}_0) \cap
L_{\infty}$.}
\\
\\
One of the main ingredients of the proof of above theorems is the
Theorem of Nakai and its consequences [1] and [17]. This motivates
characterization of solvability of finitely generated subgroups of
complex diffeomorphisms. In this direction we have the following
theorem which generalizes Theorem (C) of [5]:
\\
\\
\textbf{Theorem D.} \textit{Let $G$ be a finitely generated
subgroup of $\mathrm{Diff}(\mathbb{C}^n,0)$. If $G$ has a
polynomial growth then $G$ is solvable.}
\\

\section{Preliminaries}

Let $ \mathcal{F} $ be a (singular) foliation on $ \mathbb{C}P(2)$
and $ L \subset \mathbb{C}P(2)$ be a projective line, which is not
an algebraic solution of $\mathcal{F}$ ($L \setminus
\mathrm{Sing}(\mathcal{F})$ is not a leaf of $\mathcal{F}$). We
say that $ p \in L $ is a tangency point of $\mathcal{F}$ with
$L$, if either $p \in \mathrm{Sing} (\mathcal{F})$ or $p \notin
\mathrm{Sing}(\mathcal{F})$ and the tangent spaces of $L$ and of
the leaf of $\mathcal{F}$ through $p$, at $p$, coincide. We say
that $L$ is invariant by $\mathcal{F}$ if $\forall p \in L
\setminus \mathrm{Sing}(\mathcal{F})$, $p$ is a tangency point of
$\mathcal{F}$ with $L$. Denote by $T(\mathcal{F},L)$, the set of
tangency points of $\mathcal{F}$ with $L$.  According to [10], if
$\mathrm{Sing}(\mathcal{F})$ has codimension $\geq 2$ or
equivalently the  singularities of ${\mathcal{F}}$ are finitely
many points in $\mathbb{C}P(2)$, then there exists an open, dense
and connected subset $NI(\mathcal{F})$ of the set of lines in
$\mathbb{C}P(2)$, such that every $ L \in NI(\mathcal{F})$
satisfies the following properties:
\begin{itemize}
\item{$L$ is not invariant by $\mathcal{F}$,}
\item{$T(\mathcal{F},L)$ is an algebraic subset of $L$ defined by a polynomial of degree
     $k=k(\mathcal{F})$} in $L$ and this number is independent of $L$.
\end{itemize}
The integer $k(\mathcal{F})$ is called \textit{the degree of
foliation $\mathcal{F}$}. According to [10], a foliation of degree
$n$ in $\mathbb{C}P(2)$ can be expressed in an affine coordinate
system by a differential equation of the form
$$(P(x,y)+xg(x,y))dy-(Q(x,y)+yg(x,y))dx=0, $$ where $P$, $Q$ and
$g$ are polynomials such that:
\begin{enumerate}
\item{$P+xg$ and $Q+yg$ are relatively prime,}
\item{g is homogeneous of degree $n$,}
\item{$\mathrm{max}\{ \mathrm{deg}(P),~\mathrm{deg(Q)} \} \leq n ,$}
\item{$\mathrm{max}\{ \mathrm{deg}(P),~\mathrm{deg(Q)} \}=n$ if $g \equiv 0$.}
\end{enumerate}

Let $B_{n+1}$ be space of polynomials of degree $\leq n+1$ in two
variables. Let $V \subset B_{n+1} \times B_{n+1} $ be the subspace
of pairs of polynomials of the form $(p+xg, q+yg )$, where $P$,
$Q$ and $g$ are as in (2) and (3) above. Clearly $V$ is a vector
subspace of $B_{n+1} \times B_{n+1}$. Let $\mathbb{P}(V)$ be the
projective space of lines through $0 \in V$. Since the
differential equations $(P+xg)dy-(Q+yg)dx=0$ and $\lambda
(P+xg)dy- \lambda (Q+yg)dx=0$ define the same foliation in
$\mathbb{C}^2$, we can identify the set of all foliations of
degree $n$ in $\mathbb{C}P(2)$ with a subset $\mathcal{F}(n)
\subset \mathbb{P}(V)$. We consider $\mathcal{F}(n)$ with the
topology induced by the topology of $\mathbb{P}(V)$.
$\mathcal{F}(n)$ is called the \textit{space of foliations of
degree $n$ in $\mathbb{C}P(2)$}. We consider the following
subsets:

$\mathcal{S}(n):= \{ \mathcal{F} \in \mathcal{F}(n) |$ the
singularities of $\mathcal{F}$ are non-degenerated$\}$

$\mathcal{T}(n):=\{\mathcal{F} \in \mathcal{S}(n)  |$  any
characteristic number $\lambda $ of $\mathcal{F}$ satisfies
$\lambda \in \mathbb{C} \backslash \mathbb{Q_{+}} \}=\{\mathcal{F}
\in \mathcal{S}(n)| \mathcal{F} $ has reduced singularities$\}$

$\mathcal{A}(n):=\mathcal{T}(n)\bigcap \mathcal{X}(n)$

$\mathcal{H}(n):=\{\mathcal{F} \in \mathcal{A}(n) | $all
singularities of $\mathcal{F}$ in $L_\infty $ are hyperbolic $\}$

\begin{proposition} $\mathrm{ [10][11]}~ \mathcal{X}(n)$ is an analytic subvariety of
$\mathcal{F}(n)$ and also if $n \ge 2$ then:
\begin{enumerate}
\item {$\mathcal{T}(n)$ contains an open dense subset of $\mathcal{F}(n)$.}
\item {$\mathcal{H}(n)$ contains an open dense subset $\mathcal{M}_{1}(n)$ such
       that if $\mathcal{F} \in \mathcal{M}_1(n)$, $n\ge 2$ then:}
\begin{itemize}
\item {$ L_\infty$ is the only algebraic solution of $\mathcal{F}$}
\item {The holonomy group of the leaf $L_\infty \backslash
\mathrm{Sing}(\mathcal{F})$ is nonsolvable.}
\end{itemize}
\item {$ \mathcal{T}(n) \subset \mathcal{H}(n) \subset \mathcal{X}(n)$ are open
subsets.}
\end{enumerate}
\end{proposition}

\begin{lemma} Let $\mathcal{F} \in \mathcal{M}_1(n)$, $n \ge 2$; then each leaf
$F \ne L_\infty$ is dense in $\mathbb{C}P(2)$.
\end{lemma}

\begin{proof}
First we notice that $F$ must accumulate $L_\infty$. Since $F$ is
a non-algebraic leaf it must accumulate at some regular point $p
\in L_\infty \backslash \mathrm{Sing}(\mathcal{F})$. Choose a
small transverse disk $\Sigma \pitchfork L_\infty$ with $\Sigma
\subset V$, $V$ is a flow-box neighborhood of $p$. We consider the
holonomy group $\mathrm{Hol}(\mathcal{F},L_\infty,\Sigma)$. Then
$F$ accumulates the origin $p \in \Sigma$ and since by [15] (see
also theorem (3.1)) $G$ has dense pseudo-orbits in a neighborhood
the origin, it follows that $F$ is dense in a neighborhood of $p$
in $\Sigma$. Any other leaf $L^{\prime}$ of $\mathcal{F}$,
$L^{\prime} \ne L_\infty$ must have the same property. Using the
continuous dependence of the solutions with respect to the initial
conditions we may conclude that $F$ accumulates any point $q \in
F^{\prime}$, $\forall F^{\prime} \ne L_\infty$. Thus $F$ is dense
in $\mathbb{C}^2$ and since $L_\infty$ is $\mathcal{F}$-invariant,
$F$ is dense in $\mathbb{C}P(2)$.
\end{proof}

\begin{proposition}
Let $\{{\mathcal{F}_t}\}_ {t\in\mathbb{D}}$, $\mathcal{F}_0 =
\mathcal{F} \in \mathcal{M}_1(n)$ is an unfolding then it is
analytically equivalent to the trivial unfolding of $\mathcal{F}$
for $t \approx 0$.
\end{proposition}
\begin{proof}
Denote by $\tilde{\mathcal{F}}$ the foliation on
$\mathbb{C}P(2)\times\mathbb{D}$ such that $\forall {t}\in
\mathbb{D}$,
$\tilde{\mathcal{F}}\arrowvert_{\mathbb{C}P(2)\times\{t\}}={\mathcal{F}_t}$,

$\pi: \mathbb{C}^3\backslash\{0\}\rightarrow \mathbb{C}P(2)$ the
canonical projection and

$\Pi: (\mathbb{C}^3\backslash\{0\})\times\mathbb{D}\rightarrow
\mathbb{C}P(2)\times\mathbb{D}$ the map

$\Pi(p,t):=(\pi(p),t)$.

Denote by $\mathcal{F}^*:={\Pi}^*(\tilde{\mathcal{F}})$, pull-back
foliation on $(\mathbb{C}^3\backslash\{0\})\times\mathbb{D}$. Then
$\mathcal{F}^*$ extends to a foliation on
$\mathbb{C}^3\times\mathbb{D}$ by a Hartogs type argument.

\textit{Claim.} We may choose an integrable holomorphic 1-form
$\Omega$ which defines $\mathcal{F}^*$ on
$\mathbb{C}^3\times\mathbb{D}$ such that
\begin{eqnarray}
\Omega & = & A(x,t){d}{t}+ \sum_{i=1}^3{B}_{j}(x,t) dx_{j}
\nonumber,
\end{eqnarray}
where $B_{j}$ is a homogeneous polynomial of degree $n+1$ in $x$,
$A$ is a homogeneous polynomial of degree $n+2$ in $x$,
$\sum_{i=1}^3x_jB_j(x,t)\equiv 0$ and
$\Omega_t:=\sum_{i=1}^3B_j(x,t)dx_j$ defines
${\pi}^*({\mathcal{F}_t})$ on $\mathbb{C}^3$.

\textit{Proof of the claim.} First we remark that by triviality of
Dolbeault and Cech cohomology groups of
$\mathbb{C}^3\times\mathbb{D}$, $\mathcal{F}^*$ is given by an
integrable holomorphic 1-form, say, $\omega$ in
$\mathbb{C}^3\times\mathbb{D}$.

The restriction $\omega_t:=\omega\arrowvert_{
\mathbb{C}^3\times\{t\}}$ defines
$\mathcal{F}^*_t:={\pi}^*({\mathcal{F}_t})$ in $\mathbb{C}^3$.
Thus we may write
$\omega=\alpha(x,t)dt+\sum_{k=1}^3\beta^k(x,t)dx_k=\alpha(x,t)dt+\omega_t(x)$

Since the radial vector field $R$ is tangent to the leaves of
$\mathcal{F}^*$ we have $ \omega\circ R=0 $ so that $\omega_t
\circ R=0$, i.e. $\sum_{k=1}^3x_k\beta^k(x,t)=0$. Now we use the
Taylor expansion in the variable $x=(x_1,x_2,x_3)$ of $\omega$
around a point $(0,t)$ so that
$\omega=\sum_{j=\nu}^{+\infty}\omega_j$ where
$\omega_j(x,t):=\alpha_j(x,t)dt+\sum_{k=1}^3\beta^k_j(x,t)dx_k=\alpha_j(x,t)dt+\omega^t_j$
and $\alpha_j$, $\beta^k_j$ are holomorphic in $(x,t)$, polynomial
of degree $j$ in $x$, $\omega_\nu\equiv0$. Now the main argument
is the following:

\begin{lemma}
$\Omega=\alpha_{\nu+1}dt+\omega^t_\nu$ defines $\mathcal{F}^*$ in
$\mathbb{C}^3\times\mathbb{D}$.
\end{lemma}
\begin{proof}
Indeed, $\omega\wedge d\omega=0 \Rightarrow i_R(\omega\wedge
d\omega)=i_R(\omega).d\omega-\omega\wedge i_R(d\omega)=0$

$\omega\wedge i_R(d\omega)=0$ (since $i_R(\omega)=0$) $\Rightarrow
i_R(d\omega)=f\omega$ for some holomorphic function $f$ (Divisor
lemma of Saito). Therefore the Lie derivative of $\omega$ with
respect to $R$ is

\begin{eqnarray} L_R(\omega) & = & i_R(d\omega)+d(i_R(\omega))=f\omega .
\end{eqnarray}
On the other hand since
$\omega=\sum_{j=\nu}^{+\infty}\omega_j=\sum_{j=\nu}^{+\infty}(\alpha_j(x,t)dt+\omega^t_j)$
we obtain

\begin{eqnarray}  L_R(\omega)&=&
\sum_{j=\nu}^{+\infty}L_R(\alpha_j(x,t)dt+\omega^t_j)\nonumber\\
       &=&
\sum_{j=\nu}^{+\infty}\frac{d}{dz}[\alpha_j(e^zx,t)dt+\sum_{k=1}^3\beta^k_j(e^zx,t)e^zdx_k]\arrowvert_{z=0}
       \nonumber \\
& &(\textrm{The flow of}\quad R\quad \textrm{is}\quad
R_z(x,t)=(e^zx,t))
\nonumber\\
&=& \sum_{j=\nu}^{+\infty}[j\alpha_j(x,t)dt+(j+1)\omega^t_j]  .
\end{eqnarray}

Now we write the Taylor expansion also for $f$ in the variable
$x$. $f(x,t)=\sum_{j=0}^{+\infty}f_j(x,t)$ where $f_j(x,t)$ is
holomorphic in $(x,t)$
 homogeneous polynomial
of degree $j$ in $x$. We obtain from (1) and (2)

\begin{eqnarray}
\sum_{j=\nu}^{+\infty}j \alpha_jdt+(j+1)\omega^t_j & = &
(\sum_{k=0}^{+\infty}f_k)(\sum_{l=\nu}^{+\infty}\omega_l) \nonumber\\
& = & \sum_{j\ge \nu}(\sum_{l+k=j}f_k \omega_l) j
\alpha_jdt+(j+1)\omega^t_j \nonumber\\ &=& \sum_{l+k=j}f_k
\omega_l \nonumber\\ &=& \sum_{l+k=j}(f_k\alpha_ldt+f_k\omega^t_l)
\quad  l\ge\nu  \quad \textrm{and} \quad \forall j\ge \nu
\nonumber
\end{eqnarray}
Then
\begin{eqnarray}
j\alpha_j &=& \sum_{l+k=j}(f_k\alpha_l) \\
(j+1)\omega^t_j&=&\sum_{l+k=j}(f_k\omega^t_l) \quad \forall j\ge
\nu \quad \textrm{and} \quad l\ge\nu
\end{eqnarray}
In particular (3) and (4) imply $f_{0} \alpha_\nu=\nu \alpha_\nu $
and $f_0 \omega^t_\nu=(\nu+1)\omega^t_\nu $ then $f_0=\nu +1 ,
\alpha_\nu=0 $.

An induction argument shows that:

$j\ge \nu \Rightarrow (\alpha_{j+1}dt+ \omega^t_j) \wedge \Omega=0
$, $(\Omega:=\alpha_{\nu +1}dt+ \omega^t_{\nu})$

Finally since the degree of the foliation
$\mathcal{F}={\mathcal{F}_0}$ is $n$ we have $\nu=n+1$. This
proves the lemma (2.4).
\end{proof}

\begin{lemma} There exists a complete holomorphic vector field X on
$\mathbb{C}^3\times\mathbb{D_\epsilon}$ ,
 $\mathbb{D_\epsilon} \subset \mathbb{D}$ small subdisk, such that
 $X(x,t)=\frac{\partial}{\partial t}+ \sum_{j=1}^3
F_j(x,t)\frac{\partial}{\partial x_j}$, $\Omega \circ X=0$ and $
F_j(x,t)$ is linear on $x$.
\end{lemma}
\textit{Proof.} We may present
$\Omega=A(x,t)dt+\sum_{j=1}^{3}B_j(x,t)dx_j=A(x,t)dt+\omega_t$
where $i_R(\omega_t)=0$, $B_j$ is a homogeneous polynomial of
degree $n+1$ in $x$, $A$ is a homogeneous polynomial of degree
$n+2$ in $x$.
\\

\textit{Claim.} $\forall t\in \mathbb{D_\epsilon}$ ($\epsilon \ge
0$ small enough) we have $\mathrm{Sing}({\mathcal{F}_t})
\subset\{A( . ,t)=0\}$.

\textit{Proof of the claim.} Since $\Omega \wedge d \Omega =0$ we
have the coefficients of $dt\wedge dx_i \wedge dx_j $ equal to
zero, that is:

\begin{eqnarray}
A(\frac{\partial B_j}{\partial x_i}-\frac{\partial B_i}{\partial
x_j})+B_j\frac{\partial B_j}{\partial t}- B_i\frac{\partial
B_j}{\partial t}+ B_i\frac{\partial A}{\partial
x_j}-B_j\frac{\partial A}{\partial x_i}&=&0
\end{eqnarray}

Now given $p_0 \in \mathrm{Sing}(\mathcal{F}_{t_0})$, ($t_0
\approx 0$, so that $\mathcal{F}_{t_0} \in \mathcal{M}_1(n)$) we
have from (5) that $(B_i(p_{0},t_0)=B_i(p_0,t_0)=0):$
$A(p_0,t_0)(\frac{\partial B_j}{\partial
x_i}(p_0,t_0)-\frac{\partial B_i}{\partial x_j}(p_0,t_0))$. Since
$\mathcal{F}_{t_0}\in T(n)$ we have $\frac{\partial B_j}{\partial
x_i}(p_0,t_0) \neq \frac{\partial B_i}{\partial x_j}(p_0,t_0)(i\ne
j)$ and $A(p_0,t_0)=0$.

Using now Noether's lemma for foliations we conclude that there
exist $F_j(x,t)$ holomorphic in $(x,t)$, homogeneous polynomial of
degree $1=(n+2)-(n+1)$ in $x$, such that
$A(x,t)=\sum_{j=1}^3{F_j(x,t)B_j(x,t)}$. Now we define $X(x,t):=1
\frac{\partial}{\partial t}+\sum_{j=1}^3
F_j(x,t)\frac{\partial}{\partial x_j}$ so that $\Omega \circ
X=A-\sum_{j=1}^3 F_jB_j=0$.

In addition X is complete because each $F_j$ is of degree one in
$x$. The flow of $X$ writes $X_z(x,t)=(\Psi_z(x,t),t+z)$. Clearly
the
$\Psi_z:\mathbb{C}^3\backslash\{0\}\longrightarrow\mathbb{C}^3\backslash\{0\}$
defines an analytic equivalence between $\mathcal{F}$ and
$\mathcal{F}_z$. The proposition (1.9) is now proved.
\end{proof}

Another important remark is the following:

\begin{proposition} Let $\mathcal{F}$, $\mathcal{G}$ be foliations with hyperbolic
singularities on $\mathbb{C}P(2)$. Assume that $L_{\infty}$ is the
only algebraic leaf of $\mathcal{F}$ and that
$\mathcal{F}\arrowvert_{ \mathbb{C}^2}$ and
$\mathcal{G}\arrowvert_{\mathbb{C}^2}$ are topologically
equivalent. Then $L_{\infty}$ is also $\mathcal{G}$-invariant.
\end{proposition}

\begin{proof} Let $\phi : \mathbb{C}^2 \rightarrow \mathbb{C}^2$ be a topological
equivalence between $\mathcal{F}$ and $\mathcal{G}$ in
$\mathbb{C}^2$. We notice that given a singularity $ p \in
\mathrm{Sing}(\mathcal{F}) \cap L_{\infty}$, there exist local
coordinates $(x,y) \in U$, $x(p)=y(p)=0$, $L_{\infty} \cap
U=\{y=0\}$ such that $\mathcal{F}
\arrowvert_{U}:{x}{d}{y}-{\lambda}{y}{d}{x}=0$, $\lambda \in
\mathbb{C} \backslash \mathbb{R}$ and  $U \cap
\mathrm{Sing}(\mathcal{F})=\{p\}$. Let $U^*=U \backslash
(L_{\infty} \cap U)$, $V^*=\phi (U^*) \subset \mathbb{C}^2$,
$\Gamma :=(x=0)$, $ {\Gamma}^*:=\Gamma \cap U^* =\Gamma \backslash
\{p\}$.$\Gamma$ is the local separatrix of $\mathcal{F}$ at $p$,
transverse to $L_{\infty}$. We put ${\Gamma}^*_1=\phi({\Gamma}^*)
\subset V^*$. We remark that ${\Gamma}^*_1$ is contained in a leaf
of $\mathcal{G}$ and it is closed in $V^*$. On the other hand if
we take any local leaf $L$ of $\mathcal{F}\arrowvert_{U^*}$, $ L
\ne \Gamma $; then by the hyperbolicity of $p \in
\mathrm{Sing}(\mathcal{F})$ we have that $L$ accumulates $\Gamma$.
Thus the image $L_1=\phi(L)$ is a leaf of $\mathcal{G}
\arrowvert_{V^*}$ that accumulates ${\Gamma}^*_1 \ne L_1$.

Assume by contradiction that $L_{\infty}$ is not
$\mathcal{G}$-invariant. The curve ${\Gamma}^*_1 \subset
\mathbb{C}^2$ accumulates $L_{\infty}$. By the Flow Box Theorem, a
point of accumulation $q \in L_{\infty} \cap \bar{{\Gamma}}^*_1$
which is not a singularity of $\mathcal{G}$, must be a point near
to which the closure (in $\mathbb{C}P(2)$) $\bar{{\Gamma}}^*_1$ is
analytic.

Thus if there are no singularities of $\mathcal{G}$ in
${\bar{{\Gamma}}}^*_{1} \cap L_{\infty}$ then
${\bar{{\Gamma}}}^*_{1}$ is an algebraic $\mathcal{G}$-invariant
curve in $\mathbb{C}P(2)$. This implies that if $L_{0}$ is the
leaf of $\mathcal{F}$ on $\mathbb{C}P(2)$ that contains
${\Gamma}^*$ then $\bar{L}_0$ is an algebraic invariant curve and
$\mathcal{F}$-invariant. Since $\bar{L}_0 \ne L_{\infty} $ we have
a contradiction to our hypothesis. Therefore ${\Gamma}^*_1$ must
accumulate to some singularity $r$ of $\mathcal{G}$ in
$L_{\infty}$. Once again by the local behavior of the leaves close
to ${\Gamma}^*_1$ and due to the fact that $r$ is hyperbolic, it
follows that ${\bar{\Gamma}}^*_1$ is locally a separatrix of
$\mathcal{G}$ at $r$. Since $L_{\infty}$ is not
$\mathcal{G}$-invariant, we have two local separatrices
${\Lambda}_1$, ${\Lambda}_2$ for $\mathcal{G}$ at $r$ with
${\Lambda}_{j} \nsubseteq L_{\infty}$, $j=1,2$. Thus
${\bar{\Gamma}}^*_1$ is locally contained in ${\Lambda}_1 \cup
{\Lambda}_2$ and in particular ${\bar{\Gamma}}^*_1$ is analytic
around $r$. Since (as we have seen) ${\bar{\Gamma}}^*_1$ is also
analytic around the points $q \in \mathrm{Sing}(\mathcal{G})$, it
follows that ${\bar{\Gamma}}^*_1$ is analytic in $\mathbb{C}P(2)$
and once again it is an algebraic curve. Again we conclude that
$\Gamma$ is contained in an algebraic leaf of $\mathcal{F}$, other
than $L_{\infty}$. Contradiction!
\end{proof}
The proof given above also shows us:

\begin{proposition} Let $\mathcal{F}$, $\mathcal{G}$ be foliations on $\mathbb{C}P(2)$
both leaving invariant the line $L_{\infty}$. Let $\phi :
\mathbb{C}^2 \rightarrow \mathbb{C}^2$ be a topological invariant
equivalence for $\mathcal{F} \arrowvert_{\mathbb{C}^2}$ and
$\mathcal{F} \arrowvert_{\mathbb{C}^2}$ . Then $\phi$ takes the
separatrix set $\mathrm{S}_{\mathcal{F}}$ onto the separatix set
$\mathrm{S}_{\mathcal{G}}$.
\end{proposition}

Here $\mathrm{S}_{\mathcal{F}}$ and $\mathrm{S}_{\mathcal{G}}$ are
respectively the set of separatrices of $\mathcal{F}$ and
$\mathcal{G}$ in $\mathbb{C}^2$ that are transverse to
$L_{\infty}$ at some singular point $p \in
\mathrm{Sing}(\mathcal{F})$.

\begin{corollary} Let $\mathcal{F}_{0} \in \mathcal{H}(n)$, $n \ge 2$. Then
any $\mathbb{C}^2$-topologically trivial deformation
$\{{\mathcal{F}_t}\}_{t\in\mathbb{D}}$ of $\mathcal{F}_{0}$, is a
deformation in the class $\mathcal{H}(n)$ and it is also
$s$-trivial if we consider $\mathrm{t} \approx 0$.
\end{corollary}

\begin{proof}
First we recall that $\mathcal{H}(n)$ is open in $\mathcal{X}(n)$.
Thus it remains to use proposition (2.6) to conclude that
$\mathcal{F}_{t} \in \mathcal{H}(n)$, $\forall t \approx 0$ and
then we use proposition (2.7) to conclude that
$\{{\mathcal{F}_t}\}_ {t \approx 0}$ is $s$-trivial.
\end{proof}

\section{Fixed points and one-parameter pseudogroup}

$\mathrm{Diff}(\mathbb{C},0)$ denotes the group of germs of
complex diffeomorphisms fixing $0 \in \mathbb{C}$,
$f(z)={\lambda}{z}+ \sum_{n \ge 2 }{{a}_{n}}{{z}^{n}}$; $\lambda
\neq 0$.

Let $G \subset \mathrm{Diff}(\mathbb{C},0)$ be a finitely
generated subgroup with a set of generators ${g}_1, \cdots
,{g}_{r} \in G$ defined in a compact disk
$\bar{\mathbb{D}}_{\epsilon}$.

\begin{theorem}
 $\mathrm{[1],[15]}$ Suppose $G$ is nonsolvable. Then:
\begin{enumerate}
\item {The basin of attraction of (the pseudo-orbits of) $G$ is an open
neighborhood of the origin $\Omega$.($0 \in \Omega$)}
\item {Either $G$ has dense pseudo-orbits in some neighborhood $0 \in V \subset \Omega$ or
there exists an invariant germ of analytic curve $\Gamma $
(equivalent to $\mathrm{Im}{z}^k=0$ for some $k \in \mathbb{N}$)
where $G$ has dense pseudo-orbits and such that $G$ has also dense
pseudo-orbits in each component of $V \backslash \Gamma$.}
\item {$G$ is topologically rigid: Given another nonsolvable subgroup $G^{\prime} \subset
\mathrm{Diff}(\mathbb{C},0)$ and a topological conjugation $\phi :
\Omega \rightarrow {\Omega}^{\prime}$ between $G$ and
$G^{\prime}$, then $\phi$ is holomorphic in a neighborhood of 0.}
\item {There exists a neighborhood $0 \in W \subset V \subset \Omega$ where $G$ has a dense
set of hyperbolic fixed points.}
\end{enumerate}

\end{theorem}
\begin{remark}

 In the case $G$ is nonsolvable and contains some $f \in G$ with
${f^{\prime}(0)}^{n} \neq 1$, $\forall n \in \mathbb{Z} \backslash
\{0\}$
 (i.e., $f^{\prime}(0)= {e}^{2\pi i \lambda}, \lambda \notin \mathbb{Q}$) we
have the following:

Dense Orbits Property: There exists a neighborhood $ 0 \in V
\subset \Omega$ where the pseudo-orbits of $G$ are dense.

\end{remark}
\textbf{Holomorphic deformations in
$\mathrm{Diff}(\mathbb{C},0)$:}
\\
Let $g \in \mathrm{Diff}(\mathbb{C},0)$ defined in some open
neighborhood $0 \in \Omega$. A
 \emph{holomorphic (one-parameter) deformation} of $g$ is a map
 $G : {\mathbb{D}}_{\epsilon} \rightarrow \mathrm{Diff}(\mathbb{C},0)$,
($\epsilon > 0$)
 which verifies the four properties:
\begin{enumerate}
\item {$G(0) = g$ as germs}
\item {The Taylor expansion coefficients of $G(t)$ depend holomorphically on $t$}
\item {The radii of convergence of $G(t)$ and ${G(t)}^{-1}$ are both uniformly
 minorated by some constant $R \ge 0$ ($\forall t \in {\mathbb{D}}_{\epsilon}$)}
\item {The modules of the linear coefficient of $G(t)$ is uniformly minorated by some
constant $C \ge 0$. In particular $|({G(t)}^{-1})^{\prime}(0)|$ is
uniformly majored by $0 < t < \infty$.}
\end{enumerate}
 Given a finitely generated pseudo-group $G \subset \mathrm{Diff}(\mathbb{C},0)$
with a set of generators
 ${g}_1, \cdots {g}_{r} \in G$; a holomorphic (one parameter) deformation of
 $G$ is given by holomorphic deformation of ${g}_{j}$, $j=1, \cdots ,r$.
 We may restrict ourselves to the following situation:

 ${G}_{t}$ is an one-parameter analytic deformation of $G$ with
 $t \in \mathbb{D}$, $G_0=G$. We have
 ${g}_{1,t} \cdots {g}_{r,t}$ as a set of generators for
 ${G}_{t}$, all of them defined in a disk $\bar{\mathbb{D}}_{\delta}$ (uniformly
on
 $t$).
 We will consider dynamical and analytical properties of such deformations.
 The results we state below have their proofs reduced to the following case
which is studied in [23].

 ${g}_{1,t}(z) = {g_1}(z) + {t} {z}^{D+1}$
 where $D \in \mathbb{N}$ is fixed,

 ${g}_{2,t}(z) = {g}_{2}(z), \cdots ,{g}_{r,t}(z) = {g}_{r}(z) $.

 For such deformations we have:

 \begin{theorem}$\mathrm{[23]}$ Given a hyperbolic fixed point $p \approx 0$ for a
word
 $f={f_n} \circ {{f}_{n-1}} \circ \cdots \circ {{f}_{1}}$
 in $G$, we consider the corresponding word
 $f_t=f=f_{n,t} \circ \cdots \circ
 f_{1,t}$ in ${G}_{t}$. Then ${f}_{t}$ has a hyperbolic
 fixed point $p({t})$ given by the implicit differential equation with initial
conditions:

 \begin{eqnarray}  \frac{dp(t)}
 {{{p(t)}^{D+1}}{dt}} =
 {\frac{{f^{\prime}_t}(p(t))}{{f^{\prime}_t}(p(t))-1}}
 {{f^{\prime}_{1,t}}(p(t))}, \quad p(0) = p.\nonumber  \end{eqnarray}
\\
In particular ${p}({t})$ depends analytically in $t$ as well as
its multiplicator ${{f_t}^{\prime}}({p}({t}))$. This holds for
$|t| < \epsilon $ if $\epsilon > 0$ is small enough.

\end{theorem}

We also have:

\begin{theorem} $\mathrm{[22]}$ Let $f$ and $g$ be two non-commuting complex
diffeomorphisms defined in some neighborhood of the origin $0 \in
\mathbb{C}$, fixed by $f$ and $g$. Assume that $f'(0)= {e}^{2\pi i
\lambda}$, $g'(0)= {e}^{2\pi i \mu}$ with $\lambda ,\mu \in
\mathbb{C} \backslash \mathbb{R}$, $\mathrm{Re}\lambda
,\mathrm{Re}\mu \notin \mathbb{Q}$. Then there exist some bound
$K>0$ and some radius $r_0 >0$ such that if $r \in (0,r_0)$ and
$|t| \leq {K}{r}$ then the orbits of the pseudo-group generated by
$g$ and ${f}_{t}(z)= t + f({z}-{t})$ are dense in
$\bar{\mathbb{D}}_{r}$.
\end{theorem}

\begin{corollary}$\mathrm{[21]}$ Let $f$ and $g$ be as above. Any holomorphic
deformation of the subgroup $<f,g> \subset
\mathrm{Diff}(\mathbb{C},0)$ preserves locally at the origin the
dense orbits property.
\end{corollary}

\section{Growth of finitely generated subgroups of $ \mathrm{Diff}(\mathbb{C}^n ,0)$}
Consider a group $G$ generated by $S=\{g_1, \cdots , g_k\}$. Each
element $ g \in G $ can be represented by a word ${g_{i_1}}^{p_1}
{g_{i_2}}^{p_2} \cdots {g_{i_l}}^{p_l}$ and $\vert p_1 \vert +
\vert p_2 \vert + \cdots+ \vert p_l \vert $ is called \textit{the
length} of the word. \textit{The norm} $\Vert \gamma \Vert$
(relative to $S$) is defined as the minimal length of the words
representing $g$. Let $ B(n) $ be the set of elements $ g \in G $
with $ \Vert g \Vert \leq n $. The growth function of $g$ with
respect to the set of generators $S$ is defined as
$\gamma:=\gamma_G(n):=\vert B(n) \vert$.

We say that a function $f:\mathbb{N} \longrightarrow \mathbb{R}$
is dominated by a function $g:\mathbb{N} \longrightarrow
\mathbb{R}$, denoted by $f \preceq g$, if there is a constant $C >
0$ such that $f(n) \leq g(C.n)$ for all $n \in \mathbb{N}$. Two
functions $f,g:\mathbb{N} \longrightarrow \mathbb{R}$ are called
equivalent, denoted by $ f \sim g $, if $f \preceq g$ and $g
\preceq f$. It is known that for any two finite sets of generators
$S_1$ and $S_2$ of a group, the corresponding two growth functions
are equivalent. Note also that if $ \vert S \vert=k$, then
$\gamma(n) \leq k^n$.

Growth of group $G$ is called \textit{exponential} if $\gamma(n)
\sim e^n$. Otherwise the growth is said to be
\textit{subexponential}. Growth of group $G$ is called
\textit{polynomial} if $\gamma(n) \sim n^c$ for some $c > 0$. If
$\gamma(n) \succeq n^c$ for all $c$, the growth of $G$ is said to
be \textit{superpolynomial}. If the growth is subexponential and
superpolynomial, it is called \textit{intermediate}.

\textbf{Examples.} The finitely generated abelian groups have
polynomial growth. More precisely if $ S=\{g_1, \cdots , g_k\} $
is a minimal generating set for a free abelian group $G$ then the
growth function is
$$ \gamma(n)= \sum_{l=0}^{m}2^l {m \choose l} {{n \choose l}}. $$
Also the finitely generated nilpotent groups are of polynomial
growth [24].

The free groups with $k \geq 2$ generators have exponential
growth. Milnor and Wolf in [14] and [24] showed that a finitely
generated solvable group $G$ has exponential growth unless $G$
contains a nilpotent subgroup of finite index.

If $G$ is a finite extension of a group of polynomial growth, then
$G$ itself has polynomial growth. So we conclude if a finitely
generated group $G$ has a nilpotent subgroup of finite index then
$G$ has polynomial growth. Conversely the finitely generated
linear groups (Tits [20]), the finitely generated polycyclic
groups (Wolf [24]) and the finitely generated subgroups of
$\mathrm{Diff}(\mathbb{R}^n,0)$ (Plante-Thurston [17]) with
polynomial growth have nilpotent subgroups with finite index.
Finally in an extraordinary work [7], Gromov settled the problem
and proved if a finitely generated group $G$ has polynomial growth
then $G$ contains a nilpotent subgroup of finite index.

In [6], Grigorchuk constructed a family of groups of intermediate
growth which are the only known examples of such groups.

Denote by $\mathrm{Diff}(\mathbb{C}^n,0)$, the group of germs of
complex diffeomorphisms fixing the origin. Let $G \subset
\mathrm{Diff}(\mathbb{C}^n,0)$ be a finitely generated subgroup.
Nonsolvable finitely generated subgroups of complex diffeomorphism
in dimension $n=1$ play a fundamental role in dynamical study of
holomorphic vector fields in $(\mathbb{C}^2,0)$. In fact the
holonomy groups of irreducible components of desingularization of
the germ of a nondicritical foliation at a singularity in $
\mathbb{C}^2$  are finitely generated subgroups of
$\mathrm{Diff}(\mathbb{C}^n,0)$. Theorem of Nakai and its
consequences [1], [15] provide a new dynamical information for the
``generic'' foliation whose holonomy groups are nonsolvable as we
will see later. This motivates characterization of solvability of
finitely generated subgroups of complex diffeomorphisms. Our
objective is to prove a finitely generated subgroup of
$\mathrm{Diff}(\mathbb{C}^n,0)$ with polynomial growth is
solvable. Moreover if we know there is no finitely generated
subgroups with intermediate growth then it implies a nonsolvable
subgroup of $\mathrm{Diff}(\mathbb{C}^n,0)$ has exponential
growth.

We recall a useful lemma due to Gromov [7]:
\begin{lemma}
Let $G$ be a finitely generated (abstract) group of polynomial
growth. Then the commutator subgroup $[G,G]$ is also finitely
generated.
\end{lemma}

The following lemma is proved in [17]:
\begin{lemma} Suppose that $G$ has polynomial growth of degree $k$ and that
$$ H_0 \subset H_1 \subset \cdots  \subset H_n \subset G $$
is a finite sequence of subgroups such that for each $i$ ($ 1 \leq
i \leq n $) there is a non trivial homomorphism $f_i : H_i
\rightarrow  \mathbb{R}^{l}$ such that $ H_i \subset
\mathrm{Ker}f_i $. Then $n \leq k $.
\end{lemma}
\begin{remark}
A finitely generated solvable group of polynomial growth is
polycyclic.
\end{remark}
\textbf{Theorem D.} \textit{Let $G$ be a finitely generated
subgroup of $\mathrm{Diff}(\mathbb{C}^n,0)$. If $G$ has a
polynomial growth then $G$ is solvable.}
\begin{proof}
Put $ G_1=[G,G] $ and $ G_{i+1}=[G_i,G_i] $ for $i \in \mathbb{N}
$. By the lemma (4.2) it is enough to show that for each $G_i$
there is a non trivial homomorphism  $f_i \in \mathrm{Hom}( G_i ,
\mathbb{R}^{l} )$. Notice that $G_i \subset  \mathrm{Ker}f_{i+1}$
if there are such $f_i $'s. For simplicity denote by $H:=G_i$. By
the lemma (4.1) $H$ is finitely generated. Take a symmetric system
of generators $S:=\{h_1, \cdots , h_k\}$. We may write $ h_i(z)= z
+ \tilde{h}_i(z) $ for $i=1,\cdots ,k $  since elements of $H$ are
tangent to the identity. If $u, v \in H $ then $$ \widetilde{(v
\circ u)} (z)= \tilde{u}(z) + \tilde{v}(z) +
[\tilde{v}(u(z))-\tilde{v}(z)].$$ Write $\tilde{u}(z)=
(\tilde{u}_1(z), \cdots, \tilde{u}_{2n}(z) ) $ as a real function
and $ z=(x_1, \cdots, x_{2n}) \in \mathbb{R}^{2n}$. We have
$$(\ast) \quad \tilde{v}(u(z))=\tilde{v}(z)+ \sum^{2n}_{i=1} \tilde{u}_i(z)
\int^{1}_{0}\frac{\partial \tilde{v}}{\partial x_i} (x_1 + t
\tilde{u}_1, \cdots, x_{2n} + t \tilde{u}_{2n} ) dt.$$ We choose a
sequence $\{ z_m \}^{\infty}_{m=1} $ in $\mathbb{C}^n$ converges
to the origin such that for at least an index $ j \in \{ 1, \dots,
k\}$, all terms of the sequence $\{ \tilde{h}_j (z_m) \} $ are
different from zero. put $M_m=\mathrm{max}\{ \Vert \tilde{h}_1
(z_m) \Vert , \dots, \Vert \tilde{h}_k (z_m) \Vert \} $. $\forall
m \in \mathbb{N}, M_m > 0 $ and $ \forall i \in \{ 1, \dots, k \}
$ the sequence $\{ \tilde{h}_i (z_m)/M_m \} $ converges to, say,
$b_i$. If $h$ is an arbitrary element of $H$ such that $
\tilde{h}(z_m)/M_m $ converges to $b$, then for each generator
$h_i$ the sequence $\{ \widetilde{(h_i \circ h)}(z_m) /M_m \}$
converges to $ b+b_i $. In fact from $(\ast)$
$$\lim_{m \rightarrow \infty}M^{-1}_m(\tilde{h_i} ( h(z_m))- \tilde{h}(z_m))=0.$$
Now define $f:H \rightarrow \mathbb{R}^{2n} $ as following:

$$ h \mapsto  \lim_{m \rightarrow \infty} {\tilde{h}(z_m)}/M_m .$$
$f$ is well defined and non trivial homomorphism.

\end{proof}

\section{Proof of theorem A}

We use the terminology of [11] and some of the original ideas of
[8]. Let therefore $\{{\mathcal{F}_{t}}\}_ {t\in\mathbb{D}}$ be a
$\mathbb{C}^2$-topological trivial deformation of $\mathcal{F}_0
\in \mathcal{H}(n)$, $n\ge 2$. As we have proved in corollary
(2.8) there exists $\epsilon > 0$ such that
$\{{\mathcal{F}_{t}}\}_ {t\in {\mathbb{D}}_{\epsilon}}$ is a
$s$-trivial deformation of $\mathcal{F}_0$ in the class
$\mathcal{H}(n)$. Now we consider the continuous foliation
$\tilde{\mathcal{F}}$ on
$\mathbb{C}P(2)\times{{\mathbb{D}}_{\epsilon}}$ defined as
follows:
\begin{itemize}
\item {$\mathrm{Sing}(\tilde{\mathcal{F}})=
\bigcup_{|t|<\epsilon}{\mathrm{Sing}({\mathcal{F}}_{t})}\times{\{t\}}$
}
\item {The leaves of ${\mathcal{F}}_{t}$ are the intersections of the leaves of
$\tilde{\mathcal{F}}$ with
\\
${\mathbb{C}P(2)}\times{\{t\}}$, $\forall |t|<\epsilon$.}
\end{itemize}
Because of the topological triviality, $\tilde{\mathcal{F}}$ is a
continuous foliation on
${\mathbb{C}}^{2}\times{\mathbb{D}}_{\epsilon}$. This foliation
extends to a continuous foliation on
$\mathbb{C}P(2)\times{\mathbb{D}}_{\epsilon}$ by adding the leaf
with singularities ${L_{\infty}} \times {\mathbb{D}}_{\epsilon}$.
In order to prove that $\tilde{\mathcal{F}}$ is holomorphic we
begin by proving that it has holomorphic leaves and then it is
transversely holomorphic. This is basically done by the following
lemma:

\begin{lemma} Let ${p}_1, \cdots ,{p}_{n+1} \in L_{\infty}$ be the singularity of
$\mathcal{F}_0$ in $L_{\infty}$. Then
\begin{enumerate}
\item {There exist analytic functions $p_j(t)$, $t \in {\mathbb{D}}_{\epsilon}$
such that
\\
$\{{p}_1(t), \cdots ,{p}_{n+1}(t)\}
={\mathrm{Sing}({\mathcal{F}}_{t})} \cap {L}_{\infty}$,
${p}_{j}(0)={p}_{j}$ , $j=1, \cdots, n+1$. }

Fix  $q \in L_{\infty} \backslash
\mathrm{Sing}({\mathcal{F}}_{0})$ and take small simple loops
\\
${\alpha}_{j} \in {\pi}_1(L_{\infty} \backslash
\mathrm{Sing}({\mathcal{F}}_{0}),q)$ and a small transverse disk
$\Sigma \pitchfork L_{\infty}$. Then for $\epsilon > 0$ small we
have:

\item {The holonomy group
${G}_{t}:=\mathrm{Hol}({\mathcal{F}}_{t},L_{\infty},{\Sigma}_i)
\subset \mathrm{Diff}(\Sigma,q)$ is generated by the holonomy maps
$f_{j,t}$ associated to the loops ${\alpha}_j$ (${\alpha}_{j}$ is
also simple loop around $p_j(t)$).}

In particular we obtain

\item {${\{G_t\}}_{t \in {{\mathbb{D}}_{\epsilon}}}$ is an one-parameter
holomorphic deformation of
\\
${G}_{0}=\mathrm{Hol}({\mathcal{F}}_{0},L_{\infty},\Sigma)$. }
\item {The group ${G}_{t}$ is nonsolvable with the density orbits property, a dense
set ${\eta}_{t} \subset {\Sigma}\times {\{t\}}$ of hyperbolic
fixed points around the origin $(q,t)$. Moreover, given any ${p}_0
\in {\eta}_0$, ${p}_0={f}_0({p}_0)$, there exists an analytic
curve ${p}_{t} \in {\eta}_{t} $ such that $p(0)={p}_0$,
${f}_{t}({p}_{t})={p}_{t}$ where ${f}_t \in {G}_{t}$ is the
corresponding deformation of ${f}_0$.}
\end{enumerate}
\end{lemma}

Using above lemma we prove that $\tilde{\mathcal{F}}$ is
holomorphic close to ${L_{\infty}} \times
{\mathbb{D}}_{\epsilon}$: Given a point ${p}_0 \in {\eta}_0$ and
${f}_0 \in {G}_0$ as above, the curve ${p}(t)$ and ${f}_{t} \in
{G}_{t}$ given by (iv) above we have $\{{p}(t),|t|<\epsilon \}
\subset {\tilde{L}_{{p}_0}}\cap ({\Sigma}\times
{{\mathbb{D}}_{\epsilon}})$ where ${\tilde{L}}_{{p}_0}$ is the
$\tilde{\mathcal{F}}$-leaf through $p_0$. On the other hand
${\tilde{L}}_{p_0}$ is already holomorphic along the cuts
$\tilde{L}_{p_0} \cap ({\mathbb{C}P(2)}\times{\{t\}})$ for
$L_{{p^{\prime}}_0,t }$ for ${p}_0=({p^{\prime}}_0,0)$. This
implies that $\tilde{{L}}_{p_0}$ is analytic.

Since the curves $\{{p}(t),|t| < \epsilon \}$ with ${p}_0 \in
{\eta}_0$ are analytic and locally dense around ${{\{q\}} \times
{{\mathbb{D}}_{\epsilon}}} \subset {{\Sigma} \times
{{\mathbb{D}}_{\epsilon}}}$ it follows that any leaf $\tilde{L}$
of $\tilde{\mathcal{F}}$ is a uniform limit of holomorphic leaves
$\tilde{{L}}_{{p}_0}$ and it is therefore holomorphic. Thus
$\tilde{\mathcal{F}}$ has holomorphic leaves. We proceed to prove
that it is transversely holomorphic. This is in fact a consequence
of topological rigidity theorem [15] for nonsolvable groups of
$\mathrm{Diff}(\mathbb{C},0)$.

Fix transverse section ${\Sigma} \pitchfork {L_{\infty}}$ as
above. We may assume that $\Sigma \subset V$ where $V$ is a
flow-box neighborhood for $\mathcal{F}_0$ with $q \in V$. The
homeomorphisms ${\phi}_{t}:{\mathbb{C}}^2 \rightarrow
{\mathbb{C}}^2$ take the separatrices $S_0$ of $\mathcal{F}_0$
onto the set of separatrices $S_{t}$ of $\mathcal{F}_{t}$. Now we
use the following proposition:

\begin{proposition} Given $\mathcal{F} \in \mathcal{H}(n)$, $n \ge 2$, the set of
separatrices ${S}_{\mathcal{F}}$ of $\mathcal{F}$ is dense in
$\mathbb{C}P(2)$ and it accumulates densely a neighborhood of the
origin for any transverse disk ${\Sigma} \pitchfork L_{\infty}$,
$q \notin \mathrm{Sing}{\mathcal{F}}$.
\end{proposition}

\textit{Proof.} Indeed, given a separatrix $\Gamma \subset
{S}_{\mathcal{F}}$ the leaf $L \supset \Gamma$ is nonalgebraic for
$\mathcal{F} \in \mathcal{H}(n)$. This implies that $L \backslash
\Gamma$ accumulates $L_{\infty}$ and therefore any transverse disk
$\Sigma$ as above is cut by $L$. Now it remains to use the density
of the pseudo-orbits of $\mathrm{Hol}({\mathcal{F}},{L}_{\infty})$
stated in theorem (1.15).

Returning to our argumentation we fix any $p \in \Sigma$,
separatrix $(p_0 \in){\Gamma}_0 \subset {\mathcal{S}}_0$
 of $\mathcal{F}_0$ and denote by $P({\Gamma}_0, p)$ the local plaque of
$\mathcal{F}_0 \arrowvert_{V}$
 that is contained in ${\Gamma}_0 \cap V$ and contains the fixed point $p$.
 Put ${\Gamma}_{t}={\phi}_{t}({\Gamma}_0)$ and consider the map
 $t \mapsto p(t):=P({\Gamma}_{t},p)$. Clearly we may write
 $ p(t)={{\phi}_{t}(P({\Gamma}_0,p))} \cap {\Sigma \times {\{t\}}}$
 by choosing $\Sigma$ and $|t|$ small enough. This map $t \mapsto p(t)$ is
holomorphic
 as a consequence of proposition below:

\begin{proposition} Given any singularity $p^0_{j} \in \mathrm{Sing}{\mathcal{F}}$ there
exists a connected neighborhood $(p^0_{j} \in )\mathcal{U}_j $, a
neighborhood $\mathcal{U} \ni \mathcal{F}_0$ in $\mathcal{S}(n)$
and a holomorphic map ${\psi}_{j}: \mathcal{U} \rightarrow
\mathcal{U}_{j}$ such that $\forall \mathcal{F} \in  \mathcal{U}$,
${\psi}_{j}(\mathcal{F})=\mathrm{Sing}{\mathcal{F}} \cap
\mathcal{U}_{j}$, ${\psi}_{j}(\mathcal{F}_0)=p^0_{j}$. In
particular, if $\{{\mathcal{F}_{t}}\}_ {t \in \mathbb{D}}$ is a
deformation of $\mathcal{F}_0 \in \mathcal{H}(n)$, $n \ge 2$; then
given ${\Gamma}_0 \in S_0=S_{\mathcal{F}}$,
 ${\Sigma} \pitchfork {L_{\infty}}$, $V$ and $p \in {{\Gamma}_0 \cap \Sigma}$ as
above, there exist analytic curves
 ${p}_{j}(t)$ and $p(t)$ such that: $p_j(t)=\mathrm{Sing}{\mathcal{F}}_t  \cap
\mathcal{U}_j$, $p_j(0)=p^0_j$,
 $p(t)= P({\Gamma}_t,p(t))$, $p(0)=p$ and $p(t) \in {{\Gamma}_t \cap \Sigma }$.
\end{proposition}

Roughly speaking, the proposition says that both the singularities
and the separatrices of a foliation with nondegenerate
singularities, move analytically under analytic deformations of
the foliation.

Finally we define $h_{t}(p):=p(t)$ obtaining this way an injective
map in a dense subset of $\Sigma$ ($\mathcal{F}_0$ has dense
separatrices in $(\Sigma ,q)$), so that by the $\lambda$-lemma for
complex mapping we may extend $h_{t}$ to a map that $h_t: \Sigma
\rightarrow  \Sigma$. Moreover, it is clear that if $f_{j,t}$ is a
holonomy map as above then we have
\begin{eqnarray}
h_t(f_{j,0}(p))&=&f_{j,t}(h_{t}(p)) \nonumber
\end{eqnarray}

Because $f_0$ and $f_t$ fix the separatrices. Therefore, by
density we have $h_t \circ f_{j,0}=f_{j,t} \circ h_t$, $\forall j
\in \{1, \cdots , n+1 \}$ and the mapping $h_t$ conjugates the
holonomy groups
$G_t=\mathrm{Hol}(\mathcal{F}_t,L_{\infty},\Sigma)$ and $G_0$. By
the topological rigidity theorem $h_t$ is holomorphic which
implies that $\tilde{\mathcal{F}}$ is transversely holomorphic
close to ${L_{\infty}} \times {{\mathbb{D}}_{\epsilon}}$ [15]. The
density of $\mathcal{S}_t$, $\forall t$ assures that
$\tilde{\mathcal{F}}$ is in fact holomorphic in ${\mathbb{C} P(2)}
\times {\mathbb{D}_{\epsilon}}$.

Summarizing the discussion above we have:

\begin{proposition}Let ${\{\mathcal{F}_t\}}_{t \in \mathbb{D}}$ be a
$\mathbb{C}^2$-topologically trivial deformation of $\mathcal{F}_0
\in\mathcal{H}(n)$, $n \ge 2$. Then there exists $\epsilon > 0$
such that ${\{\mathcal{F}_t\}}_{t \in \mathbb{D}_{\epsilon}}$ is
an unfolding of $\mathcal{F}_0$ in $\mathbb{C} P(2)$.
\end{proposition}

\noindent\textit{The end of the proof of Theorem A.} The proof is
a consequence of propositions (2.3) and (5.4) above.

\section{Generalizations}

Theorem (A) may be extended to a more general class of foliations
on $\mathbb{C} P(2)$ as well as to foliations on other projective
spaces. This is the goal of this section. Before going further
into generalizations we state a kind of Noether's lemma for
foliations.

\begin{lemma}Let ${\{\mathcal{F}_t\}}_{t \in \mathbb{D}}$ be a holomorphic
unfolding of a foliation $\mathcal{F}_0$ of degree $n$ on
$\mathbb{C}P(2)$. Assume that for each singularity $p \in
{\mathrm{Sing}(\mathcal{F}_0) \cap L_{\infty}}$ the germ of
unfolding at $p$ is analytically trivial. Then there exists
$\epsilon > 0$ such that ${\{\mathcal{F}_{t}\}}_{|t| < \epsilon }$
is analytically trivial.
\end{lemma}

\begin{proof} Denote by

$\pi: \mathbb{C}^3\backslash\{0\}\rightarrow \mathbb{C}P(2)$ the
canonical projection and by

$\Pi: (\mathbb{C}^3\backslash\{0\})\times\mathbb{D}\rightarrow
\mathbb{C}P(2)\times\mathbb{D}$ the map $\Pi(p,t):=(\pi(p),t)$.

Choose a holomorphic integrable 1-form $\Omega$ which defines
$\tilde{\mathcal{F}^*}$ extension of ${\Pi}^*(\mathcal{F})$ to
$\mathbb{C}^3 \times \mathbb{D}$, so that we may choose
\begin{eqnarray}
\Omega&=&A(x,t)dt+\sum_{i=1}^3B_j(x,t)dx_j \nonumber,
\end{eqnarray}
where $A$, $B_j$ are holomorphic in $(x,t) \in \mathbb{C}^3 \times
\mathbb{D}$, homogeneous polynomial in $x$ of degree $n+2$, $n+1$;
$\sum_{i=1}^3x_jB_j =0$. The foliation ${\pi}^*({\mathcal{F}}_t)$
extends to $\mathbb{C}^3$ and this extension ${\mathcal{F}}^*_t$
is given by ${\Omega}_t =0$ for ${\Omega}_t:=
\sum_{i=1}^3B_jdx_j$.
\\

\textbf{Claim.} Given point $q \in {\mathbb{C}^3 \times
{\mathbb{D}}_{\epsilon}}$, $q \notin {{\{0\}} \times
{\mathbb{D}}}$, there exist a neighborhood $U(q)$ of $q$ in
${\mathbb{C}^3 \times {\mathbb{D}}_{\epsilon}}$ and local
holomorphic vector field $X_q \in \mathfrak{X}(U(q))$ such that $A
= {\Omega} \circ X_q$ in $U(q)$, for $\epsilon$ small enough.
\\

\textit{Proof of the claim.} If $q=(x_{1},t_{1})$ with $x_{1}
\notin {\mathrm{Sing}(\mathcal{F}_{0})}$ then $x_{1} \notin
{\mathrm{Sing}(\mathcal{F}_t)}$ for $|t|$ small enough and in
particular $x_{1} \notin {\mathrm{Sing}(\mathcal{F}_{t_{1}})}$.
Thus the existence of $X_q \in \mathfrak{X}(U(q))$ is assured in
this case. On the other hand if $x_{1} \in
{\mathrm{Sing}(\mathcal{F}_{0})}$ then we still have the existence
of $X_q \in \mathfrak{X}(U(q))$ because of the local analytical
triviality hypothesis for the unfolding at $x_{1}$.

Using the claim we obtain an open cover ${\{U_{\alpha}\}}_{\alpha
\in \mathbb{Q}}$ of $M:= {\mathbb{C}^3\backslash\{0\}} \times
{\mathbb{D}}$ with $U_{\alpha}$ connected and $X_{\alpha} \in
\mathfrak{X}(U_{\alpha})$ such that $A= {\Omega} \circ
{X_{\alpha}}$ in $U_{\alpha}$, $\forall \alpha \in \mathbb{Q}$.
Let ${U_{\alpha} \cap U_{\beta}} \neq \varnothing $ then we put
$X_{{\alpha}{\beta}}:= (X_{\alpha} - X_{\beta}) \arrowvert_{
{U_{\alpha} \cap U_{\beta}}}$ to obtain $X_{{\alpha}{\beta}} \in
\mathfrak{X}(U_{\alpha} \cap U_{\beta})$ such that ${\Omega} \circ
{X_{{\alpha}{\beta}}}=0$. Take now the rotational vector field

\begin{eqnarray}
Y&=&rot(B_{1},B_{2},B_{3})\nonumber\\
&=&(\frac{\partial B_{3}}{\partial x_{2}}-\frac{\partial
B_{2}}{\partial x_{3}}) {\frac{\partial}{\partial x_{1}}}+
(\frac{\partial B_{1}}{\partial x_{3}}-\frac{\partial
B_{3}}{\partial x_{1}}) {\frac{\partial}{\partial x_{2}}}+
(\frac{\partial B_{2}}{\partial x_{1}}-\frac{\partial
B_{1}}{\partial x_{2}}) {\frac{\partial}{\partial x_{3}}}
\nonumber.
\end{eqnarray}

$Y \in \mathfrak{X}({\mathbb{C}^3} \times {\mathbb{D}})$ and for
each $t \in {\mathbb{D}}$ we have
$i_Y(\mathrm{Vol})=d{{\Omega}_t}$ where $\mathrm{Vol}=dx_1 \wedge
dx_2 \wedge dx_3$ is the volume element of $\mathbb{C}^3$ in the
$x$-coordinates. Fixed now $q=(x_1,t_1) \notin
\mathrm{Sing}({\Omega}_{t_{1}})$ then the leaf of
$\mathcal{F}^*_{{t}_{1}}$ through $q$ is spanned by $Y(q)$ the
radial vector field $R(q)$, as a consequence of the remark above:
actually, we have ${i_R}{i_Y}(\mathrm{Vol})=i_R
 (d{\Omega}_t)=(n+1){\Omega}_t$.

Given thus $U_{{\alpha}{\beta}}:={U_{\alpha} \cap U_{\beta}} \neq
\varnothing $, since ${\Omega}_t(X_{{\alpha}{\beta}})$ we have
that $X_{{\alpha}{\beta}}$ is tangent to $\mathcal{F}^*_t$ outside
the points $(x,t) \in \mathrm{Sing}({\Omega}_t)$ so that we can
write $X_{{\alpha}{\beta}}=
{g}_{{\alpha}{\beta}}{R}+{h}_{{\alpha}{\beta}}Y$ for some
holomorphic functions $g_{{\alpha}{\beta}}$, $h_{{\alpha}{\beta}}
\in \mathcal{O}(U_{{\alpha}{\beta}} \backslash
\mathrm{Sing}({\Omega}_t))$. Since $\mathrm{Sing}({\Omega}_t)$ is
an analytic set of codimension $\ge 2$, Hartogs extension theorem
[10] implies that $g_{{\alpha}{\beta}}$, ${h}_{{\alpha}{\beta}}$
extend holomorphically to $U_{{\alpha}{\beta}}$. Now if
$U_{\alpha} \cap U_{\beta} \cap U_{\gamma} \neq \varnothing $ then
$$0=X_{{\alpha}{\beta}}+ X_{{\beta}{\gamma}} + X_{{\gamma}{\alpha}}
=(g_{{\alpha}{\beta}}+g_{{\beta}{\gamma}}+g_{{\gamma}{\alpha}})R +
(h_{{\alpha}{\beta}}+h_{{\beta}{\gamma}}+h_{{\gamma}{\alpha}}){Y}
$$ and since $R$ and  $Y$ are linearly independent outside
$\mathrm{Sing}({\Omega}_t)$ we obtain:
$g_{{\alpha}{\beta}}+g_{{\beta}{\gamma}}+g_{{\gamma}{\alpha}}=0$,
$h_{{\alpha}{\beta}}+h_{{\beta}{\gamma}}+h_{{\gamma}{\alpha}}=0$.

Thus $(g_{{\alpha}{\beta}})$, $(h_{{\alpha}{\beta}})$ are additive
cocycles in $M$ and by Cartan's theorem (for ${\mathbb{C}}^{n+1}
\backslash \{0\}$, $n \ge 2$) these cocycles are trivial, that is,
$\exists g_{\alpha}, h_{{\alpha}} \in \mathcal{O}(U_{\alpha}) $
such that if ${U_{\alpha} \cap U_{\beta}} \neq \phi $ then
$g_{{\alpha}{\beta}}=g_{\alpha} -g_{\beta}   h_{{\alpha}{\beta}}=
h_{\alpha}-h_{\beta}$ in $U_{\alpha} \cap U_{\beta}$. This gives
$X_{\alpha} - X_{\beta} = X_{{\alpha}{\beta}} =
g_{{\alpha}{\beta}}{R}+h_{{\alpha}{\beta}}{Y}
(g_{\alpha}R+h_{\alpha}Y)- (g_{\beta}R+h_{\beta}Y)$ in $U_{\alpha}
\cap U_{\beta} \neq \phi$. Thus, in $U_{\alpha} \cap U_{\beta}
\neq \phi$ we obtain $X_{\alpha}-g_{\alpha}R-h_{\alpha}Y=
X_{\beta}-g_{\beta}R-h_{\beta}Y$ and this gives a global vector
field $\tilde{X} \in \mathfrak{X}(M)$ such that $\tilde{X}
{\arrowvert}_{U_{\alpha}} := X_{\alpha}-g_{\alpha}R- h_{\alpha}Y$.
This vector field extends holomorphically to $\mathbb{C}^3 \times
{\mathbb{D}}$ and we have $( {\Omega}_t \circ \tilde{X})
\arrowvert_{U_{\alpha}}={{\Omega}_t} \circ X_{\alpha}-
g_{\alpha}{{\Omega}_t} \circ R - {h_{\alpha}}{{\Omega}_t} \circ Y
= A$ so that $ {\Omega}_t \circ \tilde{X}=A$.

It remains to prove that we may choose $\tilde{X}$ polynomial in
the variable $x$. Indeed, we write $\tilde{X}= \sum_{k=0}^{\infty}
\tilde{X}_k$ for the Taylor expansion of $\tilde{X}$ around the
origin, in the variable $x$.

Then $\tilde{X}_{k}$ is holomorphic in $(x,t)$ and homogeneous
polynomial of degree $k$ in the variable $x$. We have $A=
{\Omega}_t \circ \tilde{X}= \sum_{k=0}^{+ \infty }
{{\Omega}_t(\tilde{X}_k)}$ and since it is polynomial homogeneous
of degree $n+2$ in $x$ it follows that $k \neq 1 \Rightarrow
{\Omega}_t ( \tilde{X}_k)=0$ and ${\Omega}_t ( \tilde{X}_{1})=A$.
Since $\tilde{X}_{1}$ is linear, the flow of $\tilde{X}_{1}$ gives
an analytic trivialization for
$\{{\mathcal{F}_t}\}_{t\in{\mathbb{D}_{\epsilon}}}$.
\end{proof}

\section{Quasi-hyperbolic foliations}

Now we recall some of the features coming from [13]. A germ of
holomorphic foliation at $0 \in \mathbb{C}^{2}$, is
quasi-hyperbolic if after its reduction of singularities process
[18], we obtain an exceptional divisor that is a finite union of
invariant projective lines meeting transversely at double points
and a foliation with Saddle-type singularities: $xdy-\lambda ydx =
0$, $\lambda \in {(\mathbb{C}-\mathbb{R}) \cup \mathbb{R}_{-}}$.

In [13] we also find the notion of generic quasi-hyperbolic germs
of foliation with some dynamical restrictions on the structure of
the foliation after the reduction process.

The outstanding result is:
\begin{theorem} A topological trivial deformation of a generic quasi-hyperbolic germ
of foliation is an equisingular unfolding.

\end{theorem}

Using the concept of singular holonomy [3] we may strength this
results as follows:

\textbf{THEOREM E.}  \textit{Let $ {\{\mathcal{F}_t \}}_{t \in
\mathbb{D}}$ be a topologically trivial analytic deformation of a
germ of quasi-hyperbolic foliation $\mathcal{F}_0$ at $0 \in
\mathbb{C}^2$. We have the following possibilities: }
\begin{description}
\item[(i)] {$\mathcal{F}_0 $ admits a Liouvillian first integral and all its projective
holonomy groups are solvable,}
\item[(ii)] {$ {\{\mathcal{F}_t\}}_{|t| < \varepsilon } $ is an equisingular unfolding.}
\end{description}
\begin{proof}
Assume that all the projective singular holonomy groups of
$\mathcal{F}_0$ are solvable. In this case according to [3],
$\mathcal{F}_0$ has a Liouvillian first integral.(here we use
strongly the fact that $\mathcal{F}_0$ is quasi-hyperbolic) We may
therefore consider the case where some component of the
exceptional divisor has non solvable singular holonomy group. This
implies topological rigidity and abundance of hyperbolic fixed
points as well as the dense orbits property for this group as well
as for all the projective singular holonomy groups, which are the
main ingredients in the proof of theorem (7.1) and $
{\{\mathcal{F}_t\}}_{|t| < \varepsilon } $ is an unfolding.
\end{proof}
\textbf{Proof of Theorem B.} First we remark that by the
topological triviality on $\mathbb{C}P(2)$ we may assume that
$L_{\infty}$ is an algebraic leaf for $\mathcal{F}_t$ and that
${{\phi}_t}(L_{\infty}) = L_{\infty}$, $\forall t \in \mathbb{D}$.
In fact, we take $S_t ={\phi}_t(L_{\infty}) \subset \mathbb{C}
P(2)$. Then $S_t$ is compact $\mathcal{F}_t$-invariant and of
dimension one, so that $S_t$ is an algebraic leaf of
$\mathcal{F}_t$. By a well-known theorem of Zariski $S_t$ is
smooth. Since the self-intersection number is a topological
invariant we conclude that $S_t$ has self-intersection number one
and by Bezout's theorem $S_t$ has degree one, that is, $S_t$ is a
straight line in $\mathbb{C} P(2)$.

The problem here is that $S_t$ may do not depend analytically on
$t$. That is where we use the hypothesis that there exist at least
two reduced singularities $p_1 , p_2 \in
{\mathrm{Sing}(\mathcal{F}_0 ) \cap L_{\infty}}$. Since $p_j$ is
reduced there exists an analytic curve $p_j(t) \in
\mathrm{Sing}(\mathcal{F}_t )$ such that $p_j(t)$ is reduced
singularity of $\mathcal{F}_t$ and $ p_j(t)={\phi}_{t(p_j)}$,
$p_j(0)=p_{j}$. since the line $S_t$ contains $p_1(t) \neq p_2(t)
$ it follows that $S_t$ depends analytically on $t$ and there
exists a unique automorphism $T_t: \mathbb{C} P(2) \rightarrow
\mathbb{C} P(2)$ such that $T_t(S_t)=S_0=L_{\infty}$;
$T_t(p_j(t))=p_j$, $j=1,2$. Thus ${\psi}_t=T_t \circ {\phi}_t:
\mathbb{C} P(2) \rightarrow \mathbb{C} P(2)$ gives a topological
trivialization for the deformation ${\{{\mathcal{F}}^{1}_t\}}_{t
\in \mathbb{D}}$ of $\mathcal{F}_0$, where
${\mathcal{F}}^{1}_t:=T_t(\mathcal{F}_t)$, and $L_{\infty}$ is an
algebraic leaf of ${\mathcal{F}}^{1}_t$, $\forall t \in
\mathbb{D}$. Thus we may assume that $L_{\infty}$ is
$\mathcal{F}_t$-invariant, $\forall t \in \mathbb{D}$. Now we
proceed after performing the reduction of singularities for
$\mathcal{F}_0 \arrowvert_{L_{\infty}}$ we consider the
exceptional divisor $\mathcal{D}= \cup_{j=1}^{r}D_j$, $D_0 \cong
L_{\infty}$, $D_j \cong \mathbb{C} P(1)$, $\forall j \in \{1,
\cdots , r\}$ and observe that if the singular holonomy groups of
the components $D_j$ are all solvable then according to [3] (using
the fact that the singularities $p \in
{\mathrm{Sing}({\mathcal{F}}_{0}) \cap \mathrm{L}_{\infty}}$ are
quasi-hyperbolic) we get that $\mathcal{F}$ is a Darboux
(logarithmic) foliation. We assume therefore that some singular
holonomy group is nonsolvable, then it follows that by definition
of singular holonomy group and due to the fact that the divisor
$D$ is invariant and connected and has saddle-singularities at the
corners, we can conclude that all components of $D$ has
nonsolvable singular holonomy groups. This implies, that each germ
of ${\{\mathcal{F}_t\}}_{t \in \mathbb{D}}$ at a singular point $p
\in {\mathrm{Sing}({\mathcal{F}}_{0}) \cap L_{\infty}}$ is an
unfolding(these germs are evidently topologically trivial). Using
now arguments similar to the ones in proof of (6.1) we conclude
that ${\{\mathcal{F}_t\}}_{t \in \mathbb{D}_{\epsilon}}$ is an
unfolding for $\epsilon > 0$ small enough.

If we assume that for any singularity $p \in
{\mathrm{Sing}({\mathcal{F}}_{0}) \cap L_{\infty}}$ the germs of
unfolding is analytically trivial, then as consequence of (6.1) we
conclude that ${\{\mathcal{F}_t\}}_{t \in \mathbb{D}_{\epsilon}}$
is analytically trivial for $\epsilon > 0$ small enough. Theorem B
is now proved

\begin{remark}
 Above theorem is still true if one replace condition (3) by
the following

(3)'$ {{\phi}_t}(L_{\infty}) = L_{\infty}$, $\forall t \in
\mathbb{D}$.
\end{remark}

\section*{References}

\begin{enumerate}
\item{ Belliart M., Liousse I. and Loray F.: Sur l'existence de points fixes attractifs pour les sous-groupes de $Aut(\mathbb{C},0)$. C. R. Acad. Sci.,Paris, Ser. I, Math. 324(1997) 443-446}
\item{Camacho C., Lins Neto A. and Sad P.: Topological invariants and equidesingularization for holomorphic vector fields. J. Differential Geom. 20 (1984), no. 1, 143--174}
\item{Camacho C. and Sc$\mathrm{\acute{a}}$rdua, B.: Holomorphic foliations with Liouvillian first integrals. Ergodic Theory Dynam. Systems 21 (2001), no. 3, 717--756.}
\item{Canille Martins, Julio; Sc\'ardua, Bruno. On the growth of holomorphic projective foliations. \emph{Internat. J. Math.}  \textbf{13}  (2002),  no. 7, 695--726.}
\item{ G\'omez-Mont, X.; Ort\'iz-Bobadilla, L.; Sistemas din\'amicos holomorfos en superficies. Aportaciones Matem\'aticas: Notas de Investigaci\'on, 3. Sociedad Matem\'atica Mexicana, M\'exico, 1989.}
\item{ Grigorchuk, Rostislav I. On the Milnor problem of group growth. (Russian) \emph{Dokl. Akad. Nauk SSSR}  \textbf{271}  (1983),  no. 1, 30--33.}
\item{Gromov, Mikhael. Groups of polynomial growth and expanding maps. \emph{Inst. Hautes Études Sci. Publ. Math.} \textbf{53} (1981), 53--73.}
\item{Ilyashenko Y.: Topology of phase portraits of analytic differential equations on a complex projective plane. \emph{Trudy Sem. Petrovsk.} No. 4 (1978), 83--136}
\item{L\^e D.T.: Calcul du nombre de cycles \'evanouissaints d'une hypersurface complexe. Ann. Inst. Fourier (Grenoble) \textbf{23}  (1973) 261-270.}
\item{Lins Neto, Alcides: Algebraic solutions of polynomial differential equations and foliations in dimension two. Holomorphic dynamics (Mexico, 1986),  192--232, Lecture Notes in Math., 1345, Springer, Berlin, 1988}
\item{Lins Neto Alcides, Sad P. and Sc$\mathrm{\acute{a}}$rdua, B.: On topological rigidity of projective foliations. \emph{Bull. Soc. Math. France} \textbf{126} (1998), no. 3, 381--406.}
\item{Ma\~n\'e,, R.; Sad, P.; Sullivan, D.: On the dynamics of rational maps. Ann. Sci. \'Ec\~ole Norm. Sup. (4) 16 (1983), no. 2, 193--217.}
\item{Mattei J.F. and Salem E.: Complete systems of topological and analytical invariants for a generic foliation of $(\mathbb{C}^2,0)$, Math. Res. Lett. 4 (1) (1997), 131-141}
\item{ Milnor, John. Growth of finitely generated solvable groups. \emph{J. Differential Geometry}  \textbf{2}  1968 447--449.}
\item{Nakai, I.: Separatrices for nonsolvable dynamics on $(\mathbb{C},0)$. Ann. Inst. Fourier (Grenoble) 44 (1994), no. 2, 569--599.}
\item{Petrovski, I. and Landis, E.: On the number of limit cycles of the equation                 ${dy}/{dx}={P(x,y)}/{Q(x,y)}$ where $P$ and $Q$ are polynomials of degree 2, Amer. Math.          Society Trans., Series 2, 10 (1963) 125-176.}
\item{ Plante, J. F.; Thurston, W. P. Polynomial growth in holonomy groups of foliations.         \emph{Comment. Math. Helv.}  \textbf{51}  (1976), no. 4, 567--584. }
\item{Seidenberg, A. Reduction of singularities of the differential equation $ Ady=Bdx $.  Amer. J. Math.  90  1968 248--269.}
\item{Siu Y: Techniques of extension of analytic objects. Lecture Notes in Pure and Applied Mathematics, Vol. 8. Marcel Dekker, Inc., New York, 1974.}
\item{ Tits, J. Free subgroups in linear groups. \emph{J. Algebra}  \textbf{20} 1972 250--270. }
\item{Wirtz, B.: A property of global density for generic pseudo-groups of holomorphic diffeomorphisms of $\mathbb{C}$ with linearizable generators. Preprint 1999.}
\item{Wirtz, B.: Persistence de densite locale D'orbites de pseudo-groups de diffeomorphisms holomorphes par perturbation de points fixes. Preprint, 1999.}
\item{Wirtz, B.: Fixed points of one parameter-family of pseudo-groups of $\mathrm{Diff}(\mathbb{C},0$})
\item{ Wolf, Joseph A. Growth of finitely generated solvable groups and curvature of Riemanniann  manifolds. \emph{J. Differential Geometry}  \textbf{2}  1968 421--446. }
\end{enumerate}

\end{document}